\documentclass[12pt]{article}
\usepackage{amssymb}
\usepackage{latexsym,bm}
\usepackage{graphicx}
\usepackage{amsmath}
\usepackage{mathrsfs}

\newcommand{\bea}{\begin{eqnarray*}}
\newcommand{\eea}{\end{eqnarray*}}
\newcommand{\be}{\begin{equation}}
\newcommand{\ee}{\end{equation}}
\newcommand{\ben}{\begin{eqnarray*}}
\newcommand{\een}{\end{eqnarray*}}

\date{}
\begin{document}
\title{The minimum size of a chordal graph with given order and minimum degree
\footnote{E-mail addresses: {\tt zhan@math.ecnu.edu.cn}(X. Zhan), {\tt mathdzhang@163.com}(L. Zhang).}}
\author{\hskip -10mm Xingzhi Zhan$^a$ and Leilei Zhang$^b$\thanks{Corresponding author.}\\
{\hskip -10mm \small $^a$Department of Mathematics, East China Normal University, Shanghai 200241, China}\\
{\hskip -10mm \small $^b$School of Mathematics and Statistics, Central China Normal University, Wuhan 430079, China}
}\maketitle
\begin{abstract}
  A graph is chordal if it does not contain an induced cycle of length greater than three. We determine the minimum size of a chordal graph with given order and minimum degree.
  In doing so, we have discovered interesting properties of chordal graphs.
\end{abstract}



\section{Introduction}

We consider finite simple graphs. For terminology and notation we follow the books [3] and [5]. The {\it order} of a graph is its number of vertices, and the {\it size} is 
its number of edges. We denote by $V(G)$ and $E(G)$ the vertex set and edge set of a graph $G,$ respectively. A {\it chord} of a cycle or path $C$ is an edge not in $E(C)$ both of whose endpoints
lie on $C$. A cycle or path is called {\it chordless} if it has no chord. Chordless cycles and paths are also called {\it induced} cycles and paths, respectively.
A {\it chordal graph} is a graph in which every cycle of length greater than three has a chord. Equivalently, the graph contains no chordless cycle of length greater than three. Clearly, every induced subgraph of a chordal graph is chordal. Chordal graphs are an important family of graphs, and they are part of the larger class of perfect graphs ([1], [5, p.227]).

Extremal problems related to chordal graphs have attracted the attention of graph theorists. Blair, Heggernes, Lima and Lokshtanov [2] determined the maximum size of chordal graphs with bounded maximum degree and matching number. To study $1$-planar graphs, Zhang, Huang, Lv and Dong [6] determined the minimum size of a chordal graph with given order and connectivity.

In this paper we determine the minimum size of a chordal graph with given order and minimum degree. In doing so, we have discovered interesting new properties of chordal graphs.

We denote by $N(x)$ and $N[x]$ the neighborhood and closed neighborhood of a vertex $x$ in a graph $G$, respectively where $N[x]=N(x)\cup \{x\}.$ The degree of $x$
is denoted  by $d(x)$ or $d_G(x)$. If $H$ is a subgraph of $G,$ the {\it degree of $x$ in H,} denoted by $d_H(x),$ is defined to be $d_H(x)=|N(x)\cap V(H)|.$
For a vertex subset $S\subseteq V(G),$ we use $G[S]$ to denote the subgraph of $G$ induced by $S.$ We denote by $|G|$ and $e(G)$ the order and size of $G,$ respectively.
Thus $|G|=|V(G)|$ and $e(G)=|E(G)|.$ Denote by $\delta(G)$ the minimum degree of $G.$
For two vertex-disjoint subgraphs $H$ and $R$ of $G,$ we denote by $e(H, R)$ the number of those edges of $G$ that have one endpoint in $H$
and the other in $R.$ We use $G_1+G_2$ to denote the union of two vertex-disjoint graphs $G_1$ and $G_2,$ and use $G_1\vee G_2$ to denote the {\it join} of $G_1$ and $G_2$
which is obtained from  $G_1+G_2$ by adding edges to join every vertex in $G_1$ and every vertex in $G_2.$ $K_n$ denotes the complete graph of order $n.$
We denote by $qK_s$ the union of $q$ pair-wise vertex-disjoint complete graphs of order $s.$ An {\it (x, y)-path} is a path with endpoints $x$ and $y.$  We regard isomorphic graphs
as the same graph. Thus for two graphs $G$ and $Q,$ the notation $G=Q$ means that $G$ and $Q$ are isomorphic.

\section{Main results}

A vertex in a graph is called {\it simplicial} if its neighborhood  is a clique. We will need the following three lemmas.

{\bf Lemma 1.} (Dirac [4]) {\it Every noncomplete chordal graph has two nonadjacent simplicial vertices.}

A proof of Lemma 1 can also be found in [3, p.236]. Trivially, every vertex of a complete graph is a simplicial vertex. A vertex $x$ in a graph of order $n$ is called
{\it a dominating vertex} if $d(x)=n-1;$ otherwise it is {\it nondominating.}

{\bf Lemma 2.} {\it Let $x$ be a nondominating vertex in a chordal graph $G.$ Then there exists a vertex $y$ in $G$ such that $y$ is nonadjacent to $x$ and $G+xy$ is also a chordal graph.}

{\bf Proof.} Since the vertex $x$ is nondominating, there exists a chordless path $P_1=xy_1z_1$ of length $2.$ If $G+xz_1$ is chordal, then $y=z_1$ is the vertex we want.
Suppose $G+xz_1$ is not chordal. Then this new graph contains a chordless cycle $C$ of length at least $4.$  We must have $xz_1\in E(C).$ Denote by $P_2=xy_2z_2\ldots z_1=C-xz_1$
the path obtained from $C$ by deleting the edge $xz_1,$ which is a chordless path of length at least $3.$ Since $P_2$ is chordless and $y_1$ is adjacent to  both $x$ and $z_1,$
we have $y_1\notin V(P_2).$ Since $P_2\cup P_1$ is a cycle of length at least $5$ in $G$ and $P_2$ is chordless, $y_1$ is adjacent to every vertex of $P_2.$ In particular,
$y_1$ is adjacent to $y_2$ and $z_2.$

If $G+xz_2$ is chordal, then $y=z_2$ is what we want. If $G+xz_2$ is not chordal, then by the same argument we deduce that $G$ contains a chordless $(x,z_2)$-path $P_3=xy_3z_3\ldots z_2$ of length at least $3.$ Since $P_3$ is chordless, $y_1$ and $y_2$ are adjacent to  both $x$ and $z_2,$ we have $y_1, y_2\notin V(P_3).$ Since $P_3\cup xy_1z_2$ and
$P_3\cup xy_2z_2$ are cycles of length at least $5$ in $G$ and $P_3$ is chordless, $y_1$ and $y_2$ are adjacent to every vertex of $P_3.$ In particular, $y_1$ and $y_2$ are adjacent to $y_3$ and $z_3.$ Continuing this process we obtain a sequence $S=(P_1, P_2, \ldots)$ where for each $i,$ $P_i$ is a chordless path, for $i\ge 2,$ every $P_i=xy_iz_i\ldots z_{i-1}$ has length at least $3,$ and the set $T_i\triangleq\{y_1,y_2,\ldots, y_i\}$ is a clique of cardinality $i.$

Since $G$ is a finite graph and $|T_i|=i,$ the sequence $S$ cannot be an infinite sequence. Let $P_j$ be the last path in $S.$ Then $z_j$ is nonadjacent to $x$
and $G+xz_j$ is a chordal graph.\hfill $\Box$

{\bf Lemma 3.} {\it Let $G$ be a chordal graph of minimum degree $k$ and suppose that $p$ is a nonnegative integer with $p<k.$ Then $G$ contains a chordal spanning
subgraph of minimum degree $p.$}

{\bf Proof.} By Lemma 1, $G$ contains a simplicial vertex $x.$ We have $d(x)\ge k.$ Let $e$ be an edge incident with $x.$ We assert that $G-e$ is chordal.
To the contrary, suppose $G-e$ contains a chordless cycle $C$ of length at least $4.$ Since $C$ is also a cycle of $G$ and $G$ is chordal, $e$ must be a chord
of $C$ in $G.$ Thus, $x\in V(C).$ Let $y$ and $z$ be the two neighbors of $x$ on $C.$ Since $x$ is simplicial, $y$ and $z$ are adjacent in $G,$ and also in $G-e.$
It follows that $yz$ is a chord of $C$ in $G-e,$ contradicting the assumption that $C$ is chordless in $G-e.$

Successively deleting edges incident with $x,$ we can obtain a chordal spanning subgraph of $G$ with minimum degree $p.$ \hfill $\Box$

{\bf Theorem 4.} {\it Given positive integers $n$ and $k$ with $n\ge k+1,$ let $n=q(k+1)+r$ where $q$ and $r$ are the unique integers such that $0\leq r\leq k$.
Then the minimum size of a chordal graph of order $n$ and minimum degree $k$ is
$$
q\binom{k+1}{2}+kr-\frac{r(r-1)}{2}. \eqno (1)
$$
}

{\bf Proof.} Note that $q=\lfloor n/(k+1)\rfloor.$  Let $f(n,k)$ be the minimum size of a chordal graph of order $n$ and minimum degree $k,$ and denote the number in (1)
by $\phi(n,k).$ First we use induction on $q$ to prove the inequality
$$
f(n,k)\ge \phi(n,k). \eqno (2)
$$
Let $G$ be a chordal graph of order $n$ and minimum degree $k$. By Lemma 1, there exists a simplicial vertex $x$ in $G.$ Then $N[x]$ is a clique.
Since $\delta(G)=k,$ $|N[x]|\ge k+1.$ Let $H$ be a complete subgraph  of $G[N[x]]$ of order $k+1$ and let $R$ be the subgraph of $G$ induced by $V(G) \setminus V(H)$.
We will use the degree-sum formula $\sum_{v\in V(M)}d(v)=2 e(M)$ for any graph $M$ ([3, p.7] or [5, p.35]).

{\bf The basis step.} Suppose $q=1.$

In this case $k+1\le n\le 2k+1$ and $|R|=n-(k+1)=r.$  We have
\begin{align*}
e(G)&=e(H)+e(H,R)+e(R)\\
&=\binom{k+1}{2}+\sum_{v\in V(R)}d_{H}(v)+\frac{\sum_{v\in V(R)}d_{R}(v)}{2}\\
&=\binom{k+1}{2}+\sum_{v\in  V(R)}(d_G(v)-d_{R}(v))+\frac{\sum_{v\in  V(R)}d_{R}(v)}{2}\\
&=\binom{k+1}{2}+\sum_{v\in  V(R)}d_G(v)-\frac{\sum_{v\in  V(R)}d_{R}(v)}{2}\\
&\ge \binom{k+1}{2}+kr-\frac{r(r-1)}{2}.
\end{align*}
Note that the above argument covers also the cases $r=0,\,1.$

{\bf The induction step.} Suppose that $q\ge 2$ and inequality (2) holds for all chordal graphs of order $m$ and minimum degree $s$ such that
$q^{\prime}\triangleq\lfloor m/(s+1)\rfloor<q:$
$$
f(m,s)\ge \phi(m,s)=q^{\prime}\binom{s+1}{2}+sr^{\prime}-\frac{r^{\prime}(r^{\prime}-1)}{2} \eqno (3)
$$
where $r^{\prime}=m-q^{\prime}(s+1).$

Now $n\ge 2k+2$ and $|R|=n-(k+1)\ge k+1.$

Consider the graph $R$. Since every induced subgraph of a chordal graph is also chordal,  $R$ is chordal. We distinguish two cases.

Case 1. $\delta(R)<k.$

There exists a vertex $v$ in $R$ such that $d_R(v) < k$. Since $d_G(v)\ge k$, the vertex $v$ has $d_G(v) - d_R(v)$ neighbors in $H$. Using the fact that $|R|\ge k+1$ and applying Lemma 2, we can increase the degree of $v$ to $k$ by adding edges incident with $v$ in $R$ such that the resulting graph $R^{\prime}$ is chordal. Obviously, the number of new edges is no more than $d_G(v)-d_R(v)$. If there still exists a vertex $u$ in $R^{\prime}$ with a degree less than $k$, similarly we can increase the degree of $u$ to $k$ and keep the resulting graph chordal. The number of new edges added is no more than $d_G(u)-d_R(u)$. Continuing this process until we obtain a chordal graph $W$ with minimum degree $k$. We have
$$
e(H,R)+e(R)=\sum_{v\in V(R)}(d_G(v)-d_R(v))+e(R)\ge e(W). \eqno (4)
$$

Case 2. $\delta(R)\ge k.$

By Lemma 3, $R$ has a chordal spanning subgraph $W$ of minimum degree $k.$ Clearly $e(R)\ge e(W).$ Hence in this case, (4) holds as well.

Note that $|W|=|R|=n-(k+1)=(q-1)(k+1)+r.$ Applying the induction hypothesis (3) to the graph $W$ we have
$$
e(W)\ge \phi (n-k-1,k)=(q-1)\binom{k+1}{2}+kr-\frac{r(r-1)}{2}. \eqno (5)
$$

Using (4) and (5) we obtain
\begin{align*}
  e(G)&=e(H)+e(H,R)+e(R)\\
  &\ge e(H)+e(W)\\
  &\ge\binom{k+1}{2}+(q-1)\binom{k+1}{2}+kr-\frac{r(r-1)}{2}\\
  &=\phi (n,k).
\end{align*}
Since $G$ was arbitrarily chosen, this proves (2).

Now we show that for every pair of $(n,k),$ the lower bound $\phi(n, k)$ on $f(n,k)$ can be attained by a chordal graph $Q_{n,k}$ of order $n$ and minimum degree $k$.
Then $f(n,k)=\phi(n,k)$ and the proof will be complete. Recall that $n=q(k+1)+r$ with $0\leq r\leq k.$

If $r=0,$ let $Q_{n,k}=q K_{k+1}.$

Next suppose $r\ge 1.$

If $q=1,$ we have $k+2\le n\le 2k+1$ and $r=n-k-1.$ Let $Q_{n,k}=K_{2k-n+2}\vee (K_{n-k-1}+K_{n-k-1}).$ Note that $K_{2k-n+2}\vee K_{n-k-1}=K_{k+1}$
and $2k-n+2=k-(r-1).$

If $q=2,$ we have $2k+3\le n\le 3k+2$ and $k+2\le n^{\prime}\triangleq n-k-1\le 2k+1.$ Also $n^{\prime}-k-1\ge 1$ and $2k-n^{\prime}+2\ge 1.$
Let $e=a_1b_1$ be an edge in $Q_{n^{\prime},k}=K_{2k-n^{\prime}+2}\vee (K_{n^{\prime}-k-1}+K_{n^{\prime}-k-1})$ with $a_1$ in $K_{2k-n^{\prime}+2}$ and
$b_1$ in one $K_{n^{\prime}-k-1}.$ Add a new $K_{k+1}$ that is vertex-disjoint from $Q_{n^{\prime},k}$, and choose a vertex $c_1$ in $K_{k+1}.$
Finally deleting the edge $e$ and adding the edge $b_1c_1,$ we obtain $Q_{n,k}.$

If  $q=3,$ we have $3k+4\le n\le 4k+3.$ We will construct $Q_{n,k}$ based on $A\triangleq Q_{n-k-1,k}.$ The graph $A$ constructed above contains a cut-edge $b_1c_1.$
Add a new graph $B\triangleq K_{k+1}$ that is vertex-disjoint from $A.$ Choose two vertices $b_2,\,c_2\in V(B).$ Then deleting the edges $b_1c_1$ and $b_2c_2,$ and
adding two edges $b_2b_1$ and  $c_2c_1$, we obtain the graph $Q_{n,k}$. Note that the resulting graph $Q_{n,k}$ still has a cut-edge (the edge $c_2c_1$).

For $q\ge 4,$ we may repeat the above step from $q=2$ to $q=3$ indefinitely. This completes the proof. \hfill $\Box$

Next we consider connected chordal graphs. The case $k=1$ is trivial: The minimum size of a connected chordal graph of order $n$ and minimum degree $1$ is $n-1$
and the trees are the only extremal graphs. Thus we will assume $k\ge 2.$

\noindent{\bf Theorem 5.} {\it Given positive integers $n$ and $k$ with $n\ge k+2\ge 4,$ let $n=q(k+1)+r$ where $q$ and $r$ are the unique integers such that $0\leq r\leq k$.
Then the minimum size of a connected chordal graph of order $n$ and minimum degree $k$ is
$$
\begin{cases}
\displaystyle q\binom{k+1}{2}+kr-\frac{r(r-1)}{2} \quad {\rm if}\,\,\,r\neq 0,\\[7pt]
\displaystyle q\binom{k+1}{2}+1\quad\quad\quad\quad {\rm if}\,\,\,r=0.
\end{cases}
$$}

{\bf Proof.} Let $f(n,k)$ be defined as in the proof of Theorem 4; i.e., the number in (1), and let $g(n,k)$ be the minimum size of a connected chordal graph of order $n$
and minimum degree $k.$ Then $g(n,k)\ge f(n,k).$

We first suppose that $r\neq 0.$ Observe that the extremal graphs constructed in the proof of Theorem 4 are connected
in this case with the assumed condition $k\ge 2.$ Thus
$$
g(n,k)=f(n,k)=q\binom{k+1}{2}+kr-\frac{r(r-1)}{2}.
$$

Next suppose that $r=0.$ We have $n=q(k+1)$ and  $q\ge 2$, since $n\ge k+2.$ In this case, $f(n,k)=q\binom{k+1}{2}.$ Let $G$ be a chordal graph (not necessarily connected)
of order $n$ and minimum degree $k$ such that $e(G)=q\binom{k+1}{2}.$ Since $\delta(G)=k,$ by the degree-sum formula we have
$$
q\binom{k+1}{2}=e(G)\ge\frac{nk}{2}=q\binom{k+1}{2}.
$$
Thus $G$ must be $k$-regular. Successively using the fact that every chordal graph has a simplicial vertex (by Lemma 1), we deduce that $G=qK_{k+1}.$ Hence
$qK_{k+1}$ is the unique extremal graph for $f(n,k).$ Since $qK_{k+1}$ is disconnected, we obtain
$$
g(n,k)\ge f(n,k)+1=q\binom{k+1}{2}+1. \eqno (6)
$$

Now we recursively construct graphs $B(n,k)$ to show that the lower bound for $g(n,k)$ in (6) can be attained. Every $B(n,k)$ has a cut-edge.

If $q=2,$ we have $n=2k+2.$ The graph $B(2k+2,k)$ is obtained from $2K_{k+1}$ by adding one edge. Then $B(2k+2,k)$ is a connected chordal graph of size $2\binom{k+1}{2}+1$ with a cut-edge.

Note that $k+1\ge 3.$ Hence deleting an edge in $K_{k+1}$ yields a connected graph. Suppose $B(n,k)$ has been constructed and $a_1b_1$ is a cut-edge. Add a new graph
$F=K_{k+1}$ that is vertex-disjoint from $B(n,k).$ Let $a_2b_2$ be an edge in $F.$ Deleting the edges $a_1b_1$ and $a_2b_2$ and then adding two edges
$a_1a_2$ and $b_1b_2,$ we obtain $B(n+k+1,k).$ This completes the proof. \hfill $\Box$

{\bf Acknowledgement.} This research  was supported by the NSFC grant 12271170 and Science and Technology Commission of Shanghai Municipality
 grant 22DZ2229014.

{\bf Data availability.} No data set is used during the study.

{\bf Declarations}

{\bf Conflict of interests:} The authors have no relevant financial or non-financial interests to disclose.

\end{document}